\documentclass{article}
\usepackage{amssymb}
\usepackage{theorem}
\usepackage{amsfonts}
\title{Lagrangian Matroids:\\ Representations of Type $B_n$}
\date{17 September 2002}
\author{Richard F. Booth\thanks{Partially supported by
  The Treaty of Windsor Research Programme
of the British Council in Portugal.}
  \and Alexandre V. Borovik\thanks{Partially supported by
  The Treaty of Windsor Research Programme
of the British Council in Portugal.}
  \and Neil White\thanks{Partially supported by EPSRC grant
GR/R53593.}
        }
\theorembodyfont{\slshape}
\newtheorem{lem}{Lemma}
\newtheorem{thm}[lem]{Theorem}
\newtheorem{cor}[lem]{Corollary}

\newtheorem{ax}{Axiom}
\theorembodyfont{\upshape}
\newtheorem{defn}{Definition}
\newcommand{\qed}{\hfill $\diamond $ \\ \vskip 5pt}
\newenvironment{proof}[1][Proof]{\paragraph{#1}}{\qed}
\newcommand{\defterm}[1]{\emph{#1}}

\newcommand{\mtx}[1]{{\bf #1}}

\newcommand{\BK}{\ensuremath{\mathbb{K}}}
\newcommand{\bigp}[1]{\left( #1 \right)}

\newcommand{\bigcp}[1]{\left\{ #1 \right\}}

\newcommand{\biglset}[2]{\left\{ \left. #1 \,\right|\, #2 \right\}}

\newcommand{\B}{\ensuremath{{\cal B}}}
\newcommand{\M}{\ensuremath{{\cal M}}}

\newcommand{\sgn}{{\rm sign}\ }
\newcommand{\Pf}{{\rm Pf}}
\newcommand{\pf}{{\rm p}}

\newcommand{\tfrac}[2]{{\textstyle\frac{#1}{#2}}}
\newcommand{\bigang}[1]{\left< #1 \right>}
\newcommand{\scp}[2]{\bigang{ #1, #2}}
\newcommand\symdiff{\mathbin{\mbox{$\bigtriangleup$}}}
\newcommand{\R}{\mathbb{R}}
\newcommand{\mod}{{\rm mod}}
\begin{document}
\maketitle
\section*{Introduction}

Coxeter matroids are combinatorial objects associated with finite
Coxeter groups; they can be viewed as subsets $\cal M$ of the factor
set $W/P$ of a Coxeter group $W$ by a parabolic subgroup $P$ which
satisfy a certain maximality property with respect to a family of
shifted Bruhat orders on $W/P$.  The classical matroids of matroid
theory are exactly the Coxeter matroids for the symmetric group ${\rm
  Sym}_n$ (which is a Coxeter group of type $A_{n-1}$) and a maximal
parabolic subgroup, while the maximality property turns out to be
Gale's classical characterisation of matroids \cite{gale}.

The theory of Coxeter matroids sheds new light on the classical
matroid theory and brings into the consideration a wider class of
combinatorial objects \cite{4}.  Of this, we can specifically mention
Lagrangian matroids, which are Coxeter matroids for the
hyperoctahedral group $BC_n$ and a particular maximal parabolic
subgroup.  Lagrangian matroids are cryptomorphically equivalent to
{\em symmetric} matroids or {\em $2$-matroids} of Bouchet's papers
\cite{bou} and \cite{bou4}. They are also equivalent to {\em
  $\Delta$-matroids} \cite{bou} and to Dress and Havel's {\em
  metroids}, see \cite{bdh}.  Because of the natural embedding of
Coxeter groups $D_n < BC_n$, the {\em even} $\Delta$-matroids of
Wenzel \cite{w3} are in fact Coxeter matroids for $D_n$.

The present paper belongs to a series of publications aimed at the
development of the concept of orientation for Coxeter matroids which
would generalise the classical oriented matroids \cite{oribook}. The
concept of orientation for even $\Delta$-matroids was introduced by
Wenzel \cite{w3,wpp} and developed by Booth \cite{boo2} in a form
which better fits the general theory.  However, as we shall soon
see, this concept does not cover all natural orientation structures on
Lagrangian matroids.

When attempting to develop the theory of orientation for Coxeter
matroids other than classical (ordinary) matroids, one needs to meet
the fundamental requirement that the orientation axioms should reflect
the geometry of the appropriate flag varieties over the field $\R$ of
real numbers.  In the case of ordinary matroids of rank $k$ on $n$
elements these are the Grassmann varieties ${\mathbb G}_{n,k}$ of
$k$-dimensional subspaces in $\R^n$.  So far the two versions of
orientation of Lagrangian matroids, as developed in \cite{bbgw} and
\cite{boo2}, reflected the geometry of the flag varieties of maximal
isotropic subspaces in $C_n(\R)={\rm Sp}_{2n}(\R)$ and $D_n(\R) = {\rm
  O}_{2n,n}(\R)$. The bilinear forms on the underlying vector space
$\R^{2n}$ are, correspondingly, skew-symmetric and symmetric. Not
surprisingly, the corresponding theories of oriented Lagrangian
matroids are very different.

However, the real Lie groups $C_n(\R)={\rm Sp}_{2n}(\R)$ and
$B_{n}(\R) = {\rm O}_{2n+1,n}(\R)$ have the same Weyl group $BC_n$.
Lagrangian matroids represented in the flag variety of maximal
isotropic subspaces in the underlying vector space of $\R^{2n+1}$
${\rm O}_{2n+1,n}(\R)$ have rather natural orientation properties: as
we show in this paper, a $B_n$-represented Lagrangian matroid $\cal M$
can be obtained by gluing together a {\em Lagrangian pair\/}
\cite{Lpairs} $({\cal M}_1, {\cal M}_2)$ of $D_n$-represented
orthogonal Lagrangian matroids, and the corresponding orientation can
be very naturally described as the orientation of the {\em exploded
  sum} ${\cal M}_1 \boxplus {\cal M}_2$, which turns out to be a
$D_{n+1}$-represented Lagrangian matroid.

Hence, although $B_n$-represented Lagrangian matroids have properties
very different from that of $D_n$-represented matroids, the
corresponding orientation theory are essentially the same (up to some
non-trivial cryptomorphism).

$B_n$-representations belong to a series of representations of
Lagrangian matroids in groups ${\rm O}_{2n+m,n}$ of isometries of the
spaces $\R^{2n+m}$ endowed with non-degenerate symmetric bilinear
forms which allow maximal isotropic subspaces of dimension $n$.  At
this point we can only conjecture that these new representations are
likely to lead to the orientation theories which can be
cryptomorphically reduced to $D_n$-orientations.

\medskip

The terminology and notation follow \cite{4}.

\section{Symplectic and orthogonal matroids}

Let
$$
[n]=\{1,2,\ldots,n\} \;  \hbox{and}\; [n]^{*}=\{1^{*},2^{*},\ldots,n^{*}\}.
$$
Define the map $*:[n]\rightarrow [n]^{*}$ by $i\mapsto i^{*}$ and
the map $*:[n]^{*}\rightarrow [n]$ by $i^{*}\mapsto i $. In other
words, we are defining $i^{**}=i$. Then $*$ is an involutive
permutation of the set $[n]\cup [n]^*$.

We denote $J = [n]\cup[n]^*$.
We say that a subset $K\subset J$ is {\em admissible} if and only if
$K\cap K^{*}=\emptyset $.
If $B \subseteq J$, we set $B^+ = B \cup B^*$.

A linear ordering $\prec$ of $J$ is called a {\em
$C_n$-admissible ordering} if  $i\prec j$ implies
that\/ $j^* \prec i^*$ for all $i,j\in J$.
Equivalently, an ordering $\prec$ on $J$ is
$C_n$-admissible if and only if,
when the $2n$ elements are listed from largest to smallest,
the first $n$ elements listed form an admissible set, and the last $n$
elements listed are the stars of the first $n$ elements listed,
but are listed
in reverse order. A {\em $D_n$-admissible ordering} of $J$ is
similar to a $C_n$-admissible ordering, except that the middle two elements
(i.e., the $n$-th and $n+1$-st elements in the above listing) are now
incomparable.

Denote by $J_k$ the collection of all admissible $k$-subsets in $J$,
for some $k\leqslant n$.
If $\prec$ is $C_n$ or $D_n$-admissible ordering on $J$,
it induces the partial ordering (which we denote by the same symbol
$\prec$) on $J_k$: if $A,B \in J_k$ and
$$
A=\{a_1\prec a_2\prec \cdots \prec a_k\} \quad\hbox{ and }\quad
B=\{b_1\prec b_2\prec \cdots \prec b_k\},
$$
we set $A\prec B$ if
$$
a_1 \prec b_1, a_2 \prec b_2, \ldots, a_k \prec b_k.
$$
This partial ordering is called the Gale ordering on $J_k$ induced by $\prec$.

Now let ${\cal B} \subseteq J_k$
be a collection of admissible $k$-element subsets of the set $J$.
We say that $M=(^*,\,{\cal B})$ is a {\em
symplectic matroid\/} if it satisfies the following {\em
Maximality Property}:
\begin{quote}
{\em for every $C_n$-admissible order $\prec$ on $J$, the collection
$\cal B$ contains a unique maximal member, i.e.\ a subset $A \in
{\cal B}$ such that  $B \prec A$ {\rm (}in the Gale order
induced by $\prec${\rm )}, for all $B\in {\cal B}$}.
\end{quote}
The collection $\cal B$ is called the {\em  collection of bases} of the symplectic
matroid $M$, its elements are called {\em bases} of $M$, and the cardinality
$k$ of the bases is the {\em rank\/} of $M$. An {\em orthogonal
matroid\/} is defined similarly using $D_n$-admissible orderings.
Ordinary matroids on $[n]$ can be defined in a similar fashion, using
$A_n$-admissible orderings, which are arbitrary linear orderings on $[n]$;
indeed, this is essentially the well-known greedy algorithm of matroid theory.
A {\em Lagrangian
matroid} (resp. {\em Lagrangian orthogonal matroid\/}) is a symplectic matroid
(resp. orthogonal matroid\/) of rank $n$, the maximum possible.

One more very useful characterization of Lagrangian orthogonal matroids is the
Strong Exchange Property \cite{1}. A collection ${\cal B}\subseteq J_n$ is the
collection of bases of a Lagrangian orthogonal matroid if and only if:
\begin{quote}
For every $A,B\in{\cal B}$ and $a\in A\symdiff B$,
there exists $b\in B\smallsetminus
A$ with $b\not=a^*$, such that both $A\symdiff\{\,a,b,a^*,b^*\,\}$ and
$B\symdiff\{\,a,b,a^*,b^*\,\}$ are members of ${\cal B}$.
\end{quote}

Here, $\symdiff$ is the symmetric difference of sets.

\section{Lagrangian pairs}

Consider an admissible set of size $n-1$.  Such a set can be completed
to an admissible set of size $n$ in exactly two ways, by appending
either $i$ or $i^*$ for some $i$.  The two resulting sets are called a
\defterm{Lagrangian pair} of sets, and are characterised by the fact
that their symmetric difference is exactly $\bigcp{i,i^*}$.

Consider now two Lagrangian orthogonal matroids $\M_1$, $\M_2$ of rank
$n$ and of opposite parity.  We say that they form a Lagrangian pair
(of Lagrangian orthogonal matroids) if they satisfy:
\begin{quote}
  For every admissible ordering, the maximal bases of $\M_1$ and
  $\M_2$ under the ordering are a Lagrangian pair of sets.
\end{quote}

We say that a pair of Lagrangian subspaces of orthogonal $2n$-space
form a Lagrangian pair of subspaces if their intersection is of
dimension $n-1$.

The following result is well-known.

\begin{lem}
  \label{lem:oriflamme}
  A totally isotropic subspace of dimension $n-1$ in orthogonal
  $2n$-space is contained in exactly two Lagrangian subspaces (which
  are a Lagrangian pair).
\end{lem}

\begin{thm} {\rm \cite[Theorem~2]{Lpairs}}
  A Lagrangian pair of subspaces represent a Lagrangian pair of
  orthogonal matroids.
\end{thm}

For $\sigma \in D_n$, write $\sigma(\M)$ for the maximum basis of an
orthogonal matroid $\M$ under the ordering $\prec^\sigma$.
\begin{thm} {\rm \cite[Theorem~3]{Lpairs}}
  Given a Lagrangian pair of Lagrangian matroids, $\M_1$, $\M_2$, set
  \[
  \B = \biglset{ \sigma(\M_1)\cap\sigma(\M_2) }{ \sigma \in D_n}.
  \]
  Then $\B$ is the collection of bases of an orthogonal matroid $\M$ of rank
  $n-1$.  Furthermore, $\M_1$ and $\M_2$ are the unique Lagrangian
  orthogonal matroids obtained by completing the bases of $\M$ to
  $n$-sets of odd and even parity.
\end{thm}

\begin{thm} {\rm \cite[Theorem~7]{Lpairs}}
  Let $\B_1$, $\B_2$ be the collections of bases of a Lagrangian pair of
  Lagrangian matroids.  Then $\B=\B_1\cup\B_2$ is the collection of bases of
  a (symplectic) Lagrangian matroid.
\end{thm}

As it is shown in \cite{Lpairs}, the converse is not true.

Let $\B_1$, $\B_2$ be the collections of bases of a Lagrangian pair of
Lagrangian matroids of rank $n$.  We say that
\[
\B = \biglset{B\cup (n+1)}{\B\in \B_1}
\cup \biglset{B\cup (n+1)^*}{\B\in \B_2}
\]
is the \defterm{exploded union} of the Lagrangian pair and write
${\cal B} = \B_1 \boxplus \B_2$.

\begin{thm} {\rm \cite[Theorem~3]{Lpairs}}
  Two admissible collections of $n$-sets $\B_1$, $\B_2$ are the
  collections of bases of a Lagrangian pair of Lagrangian matroids if
  and only if their exploded union is the collection of bases of a
  Lagrangian orthogonal matroid of rank $n+1$.
\end{thm}




\section{Representations of type $D_n$ }

Concepts of representation of matroids have been introduced in two
separate, but closely related, ways. Bouchet introduces a concept of
representation by square matrices of `symmetric type' (\cite{bou3}),
whereas in \cite{5} representations are introduced in terms of
isotropic subspaces. In this paper we are concerned mainly with
representations over the real numbers.

Representable symplectic matroids arise naturally from symplectic and
orthogonal geometries, similarly to the way that classical matroids
arise from projective geometry.

\subsection{Symplectic and orthogonal representations}
Let $V$ be a vector space with basis
\[ E = \{ e_1, \ldots, e_n, e_{1^*}, , \ldots, e_{n^*} \}. \]
Let $\cdot$ be a bilinear form on $V$, with the symbol $\cdot$
often suppressed as usual, with
\begin{eqnarray*}
e_i e_{i^*} = 1 && \hbox{ for all } i \in I \\ e_i e_j = 0 && \hbox{
for all } i,j \in J \hbox{ with } i \ne j^*.
\end{eqnarray*}
\begin{defn}
  The pair $(V, \cdot)$ is called a \defterm{symplectic space} if
  $\cdot$ is antisymmetric and an \defterm{orthogonal space} if
  $\cdot$ is symmetric. If the vector space is of characteristic 2, it
  is symplectic. A subspace $U$ of $V$ is called \defterm{totally
    isotropic} if $\cdot$ restricted to $U$ is identically zero. A
  \defterm{Lagrangian subspace} is a totally isotropic subspace of
  maximal dimension (easily seen to be $n$).
\end{defn}

Choose a basis $u_1, \ldots, u_k$ of a totally isotropic subspace $U$
and represent this basis in terms of $E$, so that
\[ u_i = \sum_{j=1}^n \left( a_{i j} e_j + b_{i j} e_{j^*} \right) . \]
Now we have represented $U$ as the row space of a $k \times 2n$ matrix
$\mtx{C}= (\mtx{A}, \mtx{B})$ with columns indexed by $J$. Let \B\ be
the collection of sets of column indices corresponding to non-zero
$k\times k$ minors which are admissible; then

\begin{thm} \label{thm:repsymorthmat}
If\/ $U$ is a totally isotropic subspace of a symplectic or orthogonal
space, \B\ is the collection of bases of a symplectic or orthogonal
matroid, respectively. Note that the matroid is independent of the
choice of basis $u_1, \ldots, u_k$ of\/ $U$.
\end{thm}

This is Theorem 5 in \cite{6}; the statement for symplectic
matroids only is Theorem 2 in \cite{5}.

$C$ is called a \defterm{(symplectic/orthogonal) representation} of
$M=(J, *, \B)$, and $M$ is said to be
\defterm{(symplecticly/orthogonally) representable}. Note that
orthogonal matroids may have symplectic representations. We also note
that, when considered in matrix form, the requirement that $U$ be
totally isotropic is equivalent to the requirement that $\mtx{AB^t}$
be symmetric in the symplectic case and skew-symmetric in the
orthogonal case.

Note that we can `embed' a representation of a classical matroid as a
representation of the canonically associated Lagrangian orthogonal
matroid.  We simply make the top $k$ rows of $\mtx{A}$ (for
a matroid of rank $k$) the representation of the classical matroid,
and the remaining rows of $\mtx{A}$ zero; and the top $k$ rows of
$\mtx{B}$ zero, and the bottom $n-k$ rows an orthogonal complement of
maximal rank of $\mtx{A}$. This is clearly the required
representation, and is both a symplectic and an orthogonal
representation simultaneously.

In the case of a general, symplectically represented, symplectic
Lagrangian matroid, we assign orientations by considering essentially
signs of determinants of principal minors of the above symmetric matrices
\cite{bbgw}. Unfortunately, in skew-symmetric matrices that produces
uninteresting results, as we shall see; the correct concept is that of
the Pfaffian, which we shall define in the next section.

\subsection{Orientations}

\label{sec:ori}

Bouchet, in \cite{bou}, defines a \defterm{$\Delta $-matroid} as a
collection \B\  of
subsets of $I =[n]$, not necessarily equicardinal, satisfying the
following:
\begin{ax}[Symmetric Exchange Axiom] \label{ax:delexax}
For $A, B \in \B$ and $i \in A \symdiff B $, there exists $j \in B
\symdiff A $ such that $(A \symdiff \{i, j\}) \in \B$.
\end{ax}
It is thus immediately apparent that a classical matroid is also a
$\Delta$-matroid. Bouchet  goes on to define a \defterm{symmetric matroid}
as essentially a $\Delta$-matroid with bases extended to $n$ elements
by adding to $B \in \B$ all starred elements which do not appear,
unstarred, in $B$. Thus a symmetric matroid is a set $\B \subseteq J_n
$ satisfying:
\begin{ax} \label{ax:symexax}
For $A, B \in \B$ and $i \in A \symdiff B $, there exists $j \in B
\symdiff A $ such that $(A \symdiff \{i, j, i^* , j^* \}) \in \B$.
\end{ax}
We shall refer to these two axioms interchangeably as `the
symmetric exchange axiom' depending on the structure to which we
refer.

In this section we shall state  Wenzel's definition of (even) oriented
$\Delta$-matroids, and extend it in the obvious way to orthogonal
Lagrangian matroids. We remark parenthetically that symplectic
Lagrangian matroids (and so $\Delta$-matroids, even or otherwise) may
be oriented as described in \cite{bbgw}. We go on to discuss
representations of these objects, and prove that a representable
(classical) oriented matroid is representable as an oriented
orthogonal matroid.

\subsection{Orientation Axioms}
\label{sub:oriax}

We shall follow Wenzel in \cite{wpp} by making:
\begin{defn} \label{def:twipf}
  A map $\pf: 2^I \rightarrow \R$ is called a \defterm{twisted Pfaffian map} if
  it satisfies the following:
  \begin{enumerate}
    \item $\pf$ is not identically zero.
    \item For all $A, B \subseteq I$ with $\pf(A)\ne 0$, $\pf(B) \ne
  0$, we have $\# A = \# B \ \mod \ 2$.
    \item If $A, B \subseteq I$ and $A\symdiff B = \{i_1 < \cdots <
  i_l\}$ then we have
    \[ \sum_{j=1}^l (-1)^j \pf(A\symdiff \{i_j\}) \cdot \pf(B\symdiff
  \{i_j\}) = 0.\]
  \end{enumerate}
\end{defn}
We call two twisted Pfaffian maps \defterm{equivalent} if they differ
only by a non-zero constant scalar multiple. In fact, Wenzel makes the
definition for a `fuzzy ring' rather than for the real numbers, but we
are interested in this paper only in representations over the real
numbers. Pfaffian maps may be defined as twisted Pfaffian maps where
$\pf(\emptyset)=1$.  The Pfaffian of a square matrix of odd size is
defined to be $0$; for a $2m\times 2m$ square skew-symmetric matrix,
it is defined as follows:

\begin{defn}
Let
\[ S'_{2m} = \left\{ \sigma \in S_{2m} \mid \sigma(2k-1) = \min_{2k-1\leqslant j
\leqslant 2m} \sigma(j) \hbox{ for } 1\leqslant k \leqslant m\right\}, \]
and let $\mtx{A}$ be a skew-symmetric matrix.  Then the
\defterm{Pfaffian} of \mtx{A} is defined by
\[
\Pf\left((a_{i j})_{1 \leqslant i, j \leqslant 2m}\right) = \sum_{\sigma \in S'_{2m}} \sgn
\sigma \prod_{k=1}^m a_{\sigma(2k-1)\,\sigma(2k)}.
\]
The Pfaffian of the empty set is 1, by definition.
\end{defn}
It can be shown that the square of the Pfaffian of a (skew-symmetric)
matrix is the determinant of that matrix.

\begin{thm} \label{thm:wenzoriax}
If $\mtx{A}$ is a skew-symmetric $n \times n$ matrix, $I_1, I_2 \subseteq I$
and $I_1 \symdiff I_2 = \{ i_1,\ldots i_l \}$ with $i_j < i_{j+1}$ for $
1\leqslant j \leqslant l-1 $ then
$$
\sum_{j=1}^{l} (-1)^j \pf(I_1 \Delta \{i_j \}) \pf(I_2\Delta\{i_j
\})=0
$$
where $\pf(S) = \Pf((a_{i j})_{i, j \in S})$ for any $S \subseteq I$.
\end{thm}
This is Proposition 2.3 in \cite{w3}.

Thus a skew-symmetric matrix with real coefficients yields a Pfaffian
map, and in fact Pfaffian maps to a given ring (here, to the reals)
are in $1-1$ correspondence with skew-symmetric matrices over the same
ring (this is Theorem 2.2 in \cite{w3}). It can
be seen (see the details in \cite{boo2})
that the subsets of $I$ corresponding to non-zero
values of the twisted Pfaffian map form a $\Delta$-matroid.

We now follow \cite[Definition 2.10]{wpp} in making
\begin{defn} \label{def:oridelmat}
  An \defterm{oriented even $\Delta$-matroid} is an equivalence class of
  maps $$\pf:2^I \longrightarrow \{+1,-1,0\}$$ satisfying
  \begin{enumerate}
    \item $\pf$ is not identically zero.
    \item For all $A, B \subseteq I$ with $\pf(A)\ne 0$, $\pf(B) \ne
  0$, we have $\# A = \# B\  \mod \ 2$.
    \item If $A, B \subseteq I$ and $A\symdiff B = \{i_1 < \cdots <
  i_l\}$ and for some $w \in \{+1, -1\}$ we have
    \[ \kappa_j = w (-1)^j \pf(A\symdiff \{i_j\}) \cdot \pf(B\symdiff
  \{i_j\}) \geqslant 0\]
    for $1\leqslant j \leqslant l$, then $\kappa_j = 0$ for all $1\leqslant j\leqslant l$.
  \end{enumerate}
  We shall often speak of a map as an oriented even $\Delta$-matroid, with the
  equivalence class implicitly understood.
\end{defn}
The \defterm{bases} of the oriented even $\Delta$-matroid are those
subsets of $F\subseteq I$ for which $\pf(F) \neq 0$. We observe that
every Pfaffian map yields an oriented $\Delta$-matroid by simply
ignoring magnitudes.
\begin{lem} {\rm \cite{boo2}}
The collection of bases of an oriented $\Delta$-matroid is a
$\Delta$-matroid.
\end{lem}

 We now make the obvious definition: Take a
Lagrangian orthogonal matroid \B, with an
equivalence class of signs assigned to its bases. Two sets of signs
are said to be equivalent when they are either identical on all bases
or opposite on all bases. We express this as an equivalence class of
maps
 $$\pf:J_n\longrightarrow\{+,-,0\}$$
with $$\B = \{ A\in J_n \mid \pf (A) \neq 0\}$$ and equivalence given by $p\sim -p$.
Consider the
corresponding even $\Delta$-matroid and equivalence class of signs
$\pf'$ obtained by ignoring starred elements; that is,
$\pf' (A) = \pf(B)$, where $B\in J_n$ is the unique element with
$B\cap I=A$. Now we say that $\pf$ is an \defterm{oriented orthogonal
matroid} exactly when $\pf'$ is an oriented even $\Delta$-matroid.

\subsection{Oriented representations}
\label{sub:orirep}

\begin{thm}  {\rm \cite{wpp}} \label{thm:oridelrep}
  Given an $n\times n$ square skew-symmetric real matrix $\mtx{A}$ and
  $T\subseteq I$, define $\pf : 2^I \rightarrow \{+1,-1,0\}$ by
  setting $\pf(B)$ to be the sign of the Pfaffian of the principal
  minor indexed by $B\symdiff T$. Then $\pf$ is an oriented even
  $\Delta$-matroid, and the underlying $\Delta$-matroid is that
  represented by $\mtx{A}$ and $T$.
\end{thm}

We now move on to define a representation of an oriented orthogonal
matroid.
\begin{defn} \label{def:oriorthmat}
  Given $\mtx{C}$, an orthogonal representation of an orthogonal
  matroid $M$ over $\R$, we construct the oriented orthogonal matroid
  represented by $\mtx{C}$ as follows. Choose a basis $F$ of $M$, and
  swap columns $j$ and $j^*$ for $j\in T = F \cap I$ so that all
  columns of $F$ are in the right-hand $n$ places. Now perform row
  operations so that the right-hand $n$ columns become the identity
  matrix. Now the left-hand side, $\mtx{A'}$, is a skew-symmetric
  matrix. Since we have $\mtx{A'}$ and $T$, we have a
  representation of an oriented even $\Delta$-matroid. Unfortunately,
  this oriented even $\Delta$-matroid is dependent on the initial
  choice of $F$, although the underlying non-oriented $\Delta$-matroid
  is not, so we modify $\mtx{A'}$ as follows.

  Set \[ \varepsilon_0 = 1 \hbox{ and } \varepsilon_i = \left\{
  \begin{array}{lr} \varepsilon_{i-1} & i \notin T \\
  -\varepsilon_{i-1} & i \in T. \end{array} \right. \] for $i>0$. Then
  set $a_{ij} = \varepsilon_i \varepsilon_j a'_{ij}$.
  $\mtx{A}=(a_{ij})$ is again skew-symmetric, with rows and columns
  indexed by $I$, and we assign to the basis $B$ the sign of the
  Pfaffian of the principal minor of $\mtx{A}$ indexed by $(B\symdiff
  F)\cap I$. If we consider instead that we have permuted column
  labels with columns, then the indices giving rise to this Pfaffian
  are those of the columns of $\mtx A$ labelled by elements of $B$.
  Note that this corresponds to the oriented even $\Delta$-matroid
  represented by $\mtx{A}, T$.
\end{defn}

Notice that this definition may be rather simply stated as follows:
\begin{itemize}
\item Standard row operations are permitted.
\item Swapping columns $i$ and $i^*$, and the associated column
  labels, is permitted after multiplying all columns $i,\ldots,n$ and
  $i^*, \ldots, n^*$ by $-1$.
\item If the right-hand $n$ columns of the representation form an
  identity matrix, write $T$ for the admissible $n$-set of their
  column indeces.  Now the left-hand $n$ columns form a skew-symmetric
  matrix $\mtx A$, and we assign signs as in the underlying
  $\Delta$-matroid represented by $\mtx A$ and $T$.
\end{itemize}

\begin{thm}  {\rm \cite{boo2}} \label{thm:oriorthmat}
  The above procedure obtains an oriented orthogonal matroid, which is
  independent of choice of $F$.
\end{thm}

\section{Representations of type $B_n$}

As usual, we write $J=[n]^+= [n]\cup[n]^*$.  We shall also use the
index set $K=\bigcp{2n+1,\ldots,2n+m}$.

We begin with a \emph{standard orthogonal space} $V^{2n+m}$, which is
a vector space $V$ over $\BK$ with basis
\[
E=\bigcp{e_1,e_2,\ldots,e_n,
  e_{1^*},e_{2^*},\ldots,e_{n^*},
  f_{2n+1},\ldots, f_{2n+m}
}
\]
and which is endowed with a symmetric bilinear form $\scp{\ }{\ }$
(which we shall call the \emph{scalar product}) such that
$\scp{e_i}{e_j}=0$ for all $i,j\in J,i\not=j^*$, $\scp{e_i}{f_k} =0$
for all $i \in J$ and $k \in K$ whereas $\scp{e_i}{e_{i^*}}=1 $ for
$i\in [n]$ and $\scp{f_k}{f_k}=1$ for $k \in K$.

A \emph{totally isotropic subspace} of $V$ is a subspace $U$ such that
$\scp uv=0$ for all $u,v\in U$. Let $U$ be a totally isotropic
subspace of $V$ of dimension $l$.

The following is one of the  standard facts on symmetric bilinear forms.

\begin{lem}
  Assume that $\BK$ is either formally real or $m=1$. Then $l\leqslant n$.
\end{lem}

Now choose a basis $\{u_1,u_2,\ldots,u_l\}$ of $U$, and expand
each of these vectors in terms of the basis $E$:
\[
u_i=\sum_{j=1}^na_{i,j}e_j+\sum_{j=1}^n
b_{i,j}e_{j^*}+\sum_{k=1}^m c_{i,k}f_k.
\]
Thus we have represented the totally isotropic subspace $U$ as
the row-space of a $l\times (2n+m)$ matrix $(A,B,C)$,
$A=(a_{i,j})$, $B=(b_{i,j})$, $C= (c_{i,k})$ with the columns
indexed by $J\cup K$, specifically, the columns of $A$ by $[n]$,
those of $B$ by $[n]^*$, and those of $C$ by $K$.

A direct computation proves
\begin{lem}
  A subspace $U$ of the standard orthogonal space $V$ is totally
  isotropic if and only if $U$ is represented by a matrix $(A,B,C)$
  with
  \begin{equation}
    AB^t+BA^t+CC^t=0. \label{eq:identity}
  \end{equation}
\end{lem}

Now, given a $k\times (2n+m)$ matrix $D=(A,B,C)$ with columns indexed
by $J\cup K$, let us define a family $\B \subseteq J_k$ by
saying $X\in \B$ if $X$ is an admissible $k$-set and the
$k\times k$ minor formed by taking the columns of $D$ indexed by
elements of $X$ is non-zero.

\begin{lem}
  \label{roweq}
  Let $D=(A,B,C)$ be a matrix defining a family ${\cal B}$, and let
  $D^\prime$ be a matrix which is row-equivalent to $D$.  Then
  $D^\prime$ defines the same family ${\cal B}$. If $D$ satisfies the
  identity {\rm (\ref{eq:identity})}, then $D'$ satisfies the same
  identity.
\end{lem}

\begin{proof}
  Elementary row operations do not change the dependencies among
  columns of $D$, hence they do not change which $k\times k$ minors
  are non-zero. Furthermore, they do not change the row-space of $D$,
  hence the total isotropy  of the
  corresponding subspace $U$ of $V$, and therefore the identity
  {\rm(\ref{eq:identity})}.
\end{proof}

\begin{thm}
\label{rep:Bn} Assume that either\/ $\BK$ is formally real, or\/ $m=1$.
If\/ $U$ is totally isotropic, then $\B$ is the collection
of bases of a symplectic matroid.
\end{thm}

\begin{proof}
  Let $D$ be the matrix corresponding to $U$, and let $\prec$ be an
  admissible order on $J$. We must show that ${\cal B}$ has a unique
  maximal member.  Let ${\cal A}$ be the collection of all
  $k$-elements subsets of $J$ (admissible or not) such that the
  corresponding $k\times k$ minor of $D$ is non-zero.  In fact we will
  show that ${\cal A}$ has a unique maximal member, and that this
  member is also in ${\cal B}$, and is therefore clearly the unique
  maximal member of ${\cal B}$, since ${\cal B}$ is a subcollection of
  ${\cal A}$.

  Let us reorder the columns of $D$ according to the order $\prec$ on
  their indices, starting with the largest index.  Let $E$ be the
  row-echelon form of that matrix.  Let $A$ be the set of indices of
  the $k$ pivot columns of $E$. Clearly $A$ is the unique maximal
  member of ${\cal A}$.  Suppose that $A$ is not admissible. Thus we
  may assume that $j,j^*\in A$, and consider the two rows $r_t$ and
  $r_q$ of $E$ in which the non-zero entries of the columns indexed by
  $j$ and $j^*$ occur:
\[
\bordermatrix{&a&b&\ldots&j&\ldots&t&t^*&\ldots&j^*&\ldots&b^*&a^*
&2n+1 & \ldots & 2n+m \cr
&0&0&\ldots&1&\ldots&*&*&\ldots&0&\ldots&*&* &* & \ldots & * \cr
&0&0&\ldots&0&\ldots&0&0&\ldots&1&\ldots&*&* &* & \ldots & *
\cr}.
\]
Since $D$ is isotropic, the scalar product $\scp{r_q}{r_q}$ of the
$q$th row with itself is zero. On the other hand, tracing the way in
which the scalar product is calculated, we immediately see that it
equals
\[
\scp{r_q}{r_q} = c_{q,2n+1}^2+\cdots + c_{q,2n+m}^2.
\]
Since either $\BK$ is formally real or $m=1$, we conclude that all
$c_{q,k}=0$, and our two rows look like this:
\[
\bordermatrix{&a&b&\ldots&j&\ldots&t&t^*&\ldots&j^*&\ldots&b^*&a^*
&2n+1 & \ldots & 2n+m \cr
&0&0&\ldots&1&\ldots&*&*&\ldots&0&\ldots&*&* &* & \ldots & * \cr
&0&0&\ldots&0&\ldots&0&0&\ldots&1&\ldots&*&* &0 & \ldots & 0 \cr}.
\]
But now it is easy to see that $\scp{r_q}{r_q} =1$, which contradicts
our assumption that $U$ is isotropic.
\end{proof}

A symplectic matroid $\cal B$ which arises from a matrix $(A,B,C)$,
with $AB^t+BA^t+CC^t=0$, is called a
\defterm{${\rm O}(2n+m,n)$-representable symplectic matroid}, and $(A,B,C)$
(with its columns indexed by $J\cup K$) is a \defterm{representation}
or \defterm{coordinatisation} of it (over the field $\BK$). If $m=1$,
we shall call $\cal B$ a \defterm{$B_n$-representable symplectic
  matroid}.

\begin{thm}
  The union of the Lagrangian pair of Lagrangian matroids represented
  by a Lagrangian pair of subspaces is a $B_n$-represented Lagrangian
  matroid.  Furthermore, every $B_n$-represented Lagrangian matroid
  either arises in this way, or is itself a represented Lagrangian
  orthogonal matroid.
\end{thm}
\begin{proof}
  Consider a pair of Lagrangian subspaces.  We may represent them as
  \[
  \bigp{
    \begin{array}{c}
      A \\ x
    \end{array}
  }\quad\hbox{and}\quad
  \bigp{
    \begin{array}{c}
      A \\ y
    \end{array}
  }
  \]
  where $A$ is an $(n-1)\times 2n$ matrix and $x$ and $y$ are $1\times 2n$
  row vectors.  Now, every row of $A$ is orthogonal to itself, every
  other row of $A$ and each of $x$ and $y$.  Choose some
  $c,\alpha,\beta \in \BK$ such that $\scp{\alpha x}{\beta y} =
  -\tfrac 12 c^2$.  Now the matrix
  \[
  \bigp{
    \begin{array}{cc}
      A & 0 \\ \alpha x+\beta y & c
    \end{array}
  }
  \]
  is a totally isotropic subspace of $V$.  Consider an $n$-subset of
  the column indices $1,\ldots, 2n$, and the $n\times n$ minor
  corresponding to these columns.  Its determinant is the sum of the
  determinants of the corresponding minors in the Lagrangian pair, not
  both of which are non-zero (since they represent matroids of
  opposite parity).  Thus, the collection of bases produced is exactly
  the union of the collections produced from the Lagrangian pair.

  For the converse, suppose we are given a $B_n$-representation.  If
  the $(2n+1)$th column is empty, then removing it we have an orthogonal
  representation of a Lagrangian orthogonal matroid, and are done.
  Otherwise, we perform elementary row operations to obtain the form
  \[
  \bigp{
    \begin{array}{cc}
      A & 0 \\ z & 1
    \end{array}
  }.
  \]
  Since $A$ spans a totally isotropic subspace $U$ of dimension $n-1$,
  it is contained in a unique pair of Lagrangian subspaces, which form
  a Lagrangian pair, by Lemma~\ref{lem:oriflamme}; again, we shall
  write $x,y$ for row vectors completing $A$ to matrices spanning each
  of these two spaces.  Since $U$ is $n-1$-dimensional, its annulator is of dimension $n+1$, and so is generated by the rows of
  $A$ and the vectors $x$ and $y$.  Since $z$ is orthogonal to every
  row of $A$ and not contained in $U$, it can be expressed in the form
  $\alpha x+\beta y$ (after some row operation), and so the Lagrangian
  pair of subspaces represent a Lagrangian pair of matroids whose
  union is our $B_n$-represented Lagrangian matroid by the same
  argument as above.
\end{proof}

\section{Orientations of type $B_n$}

Consider a $B_n$-represented matroid \M{} with basis collection \B{}
which has the basis $[n]^*$.  (If not, we can simply
swap columns $i$ and $i^*$ for appropriate choices of $i$ to obtain
such a basis.)  After row operations, its representation has the form
\[
\bigp{
  \begin{array}{ccc}
    A & I_n & c
  \end{array}
},
\]
where $A$ is an $n\times n$ matrix and $c$ a $n\times 1$ column
vector.  Since the row space of the matrix is a totally isotropic
subspace, we obtain
\[
A + A^T + c c^T = 0, \quad\hbox{ yielding }\quad S + S^T = 0
\hbox{ where } S = A - \tfrac 12 c c^t.
\]
Thus bases of \M{} correspond to non-zero determinants of diagonal
minors of the matrix $A = S+\tfrac 12 cc^t$.
\begin{thm}
  \label{thm:rep-by-dn-matrix}
  The determinant of the minor of $A=S+\tfrac 12 cc^t$ indexed by $I$
  is exactly the determinant of the minor of the skew-symmetric matrix
  \[
  \bigp{
    \begin{array}{cc}
      S & c \\
      -c^T & 0
    \end{array}
  }
  \]
  indexed by $I$ (if it is of even cardinality) or $I\cup(n+1)$ (if $I$ is of
  odd cardinality).
\end{thm}
\begin{proof}
  For notational convenience, we shall notate the appropriate minors
  of $A$, $S$ and $c$ as though they were the full matrices.

  Clearly,
  \begin{eqnarray*}
    \det A &=&
    \det\bigp{
      \begin{array}{cc}
        S+cc^t & \\
        0 & 1
      \end{array}}\\
    &=&
    \det\bigp{
      \begin{array}{cc}
        S+cc^t & c\\
        0 & 1
      \end{array}}
    \quad\hbox{(by row operations)}\\
    &=&
    \det\bigp{
      \begin{array}{cc}
        S & c\\
        -c^T & 1
      \end{array}}
    \quad\hbox{(by column operations)}\\
    &=&
    \det\bigp{
      \begin{array}{cc}
        S & c\\
        -c^T & 0
      \end{array}}
    +
    \det\bigp{
      \begin{array}{cc}
        S & c\\
        0 & 1
      \end{array}}
    \quad\hbox{(by standard identities)}\\
    &=&
    \det\bigp{
      \begin{array}{cc}
        S & c\\
        -c^T & 0
      \end{array}}
    +
    \det\bigp{S}
    \quad\hbox{(by expansion).}\\
  \end{eqnarray*}
  Now, both these determinants are skew-symmetric, and so are non-zero
  only if of even cardinality.  This completes the proof.
\end{proof}
\begin{cor}
  Append $n+1$ and $(n+1)^*$ to the bases of \M{} so that all
  resulting sets have an even number of unstarred elements.  Since
  this is one of the two possibilities for the exploded sum of the
  Lagrangian pair of matroids corresponding to \M{}, this produces the
  collection of bases of an orthogonal Lagrangian matroid which contains
  the basis $[n+1]^*$.  Then this orthogonal Lagrangian matroid is
  represented by
  \[
  \bigp{
    \begin{array}{cccc}
      S & c & I_n & 0\\
      -c^T & 0 & 0 & 1
    \end{array}
  }.
  \]
\end{cor}
Since, as with orthogonal matroids, our determinants arise from
skey-symmetric matrices, the signs of these matrices cannot possibly
be interesting; they are determined only be rank.  We again turn to
the Pfaffian.
\begin{defn}
  Given a $B_n$-representation of a Lagrangian matroid \M{}, we define
  the signs of the bases of \M{} according to the following procedure:
  \begin{itemize}
  \item If \M{} is represented by a matrix of the form $(A,I_n,c)$,
    with the columns of $A$ indexed by $D^*$ and those of $I$ by $D$,
    for some $D\in J_n$, then
    \[
    \bigp{
      \begin{array}{cccc}
        S & c & I_n & 0\\
        -c^T & 0 & 0 & 1
      \end{array}
    },
    \]
    represents an orthogonal Lagrangian matroid which is an explosion
    of \M, with columns labelled
    $D^*, w^*, D, w$, where $w\in\bigcp{n+1,(n+1)^*}$ is chosen so that
    $(D\cup w)\cap [n+1]$ has even cardinality.  Now the signs of the bases
    of \M{} are the signs of the corresponding bases of this new
    matroid.
  \item Given a representation of \M{}, we may swap columns and column
    labels $i$ and $i^*$, provided we multiply columns with labels $$i,
    i^*, i+1, (i+1)^*, \ldots, n, n^*$$ by $-1$.
  \item We can perform any standard row operations.
  \end{itemize}
\end{defn}

\begin{thm}
  The signs given by the procedure above are independent of the choice
  of\/ $D$, up to global sign change.
\end{thm}

\begin{proof}
  Since the rules for swapping columns are the same as those in the
  resulting $D_n$-orientation from definition \ref{def:oriorthmat}, it
  is enough to prove that swapping only columns $n$ and $n^*$ gives
  the correct relative signs.  So, we assume that we have a
  $B_n$-representation in \emph{canonical form}:
  \[
  \bigp{
    \begin{array}{cc|cc|c}
      S -\tfrac{1}{2}x^2 bb^T & a -\tfrac{1}{2}x^2b &I_{n-1}&0&xb\\
      -a^T -\tfrac{1}{2}x^2 b^T & -\tfrac{1}{2}x^2 &0&1& x
    \end{array}
  }.
  \]
  Here $S$ is an $(n-1)$-square skew symmetric matrix, $a$ and $b$ are
  column vectors of dimension $n-1$, and $x$ is a constant; any
  $B_n$-representation in which both $[n]^*$ and $[n-1]^*\cup n$ are
  bases can be written in this way by performing row operations to
  obtain an identity matrix in the (necessarily independent) set of
  columns $[n]^*$.  The rest of the structure shown follows from
  considering the row-orthogonality of the matrix.

  From the definition, our signs come from the $(n+1)$-square skew
  symmetric matrix
  \[
  \bigp{
    \begin{array}{ccc}
      S & a & xb \\
      -a^T & 0 & x \\
      -xb^T & -x & 0
    \end{array}
  }.
  \]

  Now we consider swapping columns $n$ and $n^*$ before expanding the
  matrix.  Our representation becomes
  \[
  \bigp{
    \begin{array}{cc|cc|c}
      S -\tfrac{1}{2}x^2 bb^T & 0  &I_{n-1}&-a -\tfrac{1}{2}x^2b & xb\\
      -a^T -\tfrac{1}{2}x^2 b^T & -1 &0& \tfrac{1}{2}x^2 & x
    \end{array}
  }.
  \]
  Inverting this right-hand-side, we obtain
  \[
  \bigp{
    \begin{array}{cc|cc|c}
      T -\tfrac{1}{2}y^2 bb^T & b -\tfrac{1}{2}y^2a &I_{n-1}&0&ya\\
      -b^T -\tfrac{1}{2}y^2 a^T & -\tfrac{1}{2}y^2 &0&1& y
    \end{array}
  },
  \]
  where $y = 2/x$ and $T = S+ba^T-ab^T$, another skew-symmetric
  matrix.  Thus the signs now come from
  \[
  \bigp{
    \begin{array}{ccc}
      T & b & ya \\
      -b^T & 0 & y \\
      -ya^T & -y & 0
    \end{array}
  }.
  \]
  Here the columns are indexed by $1,\ldots, n-1, n^*, (n+1)^*$.
  Consider the sign of a basis, in each of these two $n+1$-square
  matrices.  We take four cases:
  \begin{enumerate}
  \item The basis contains neither $n$ nor $n+1$.  Thus the signs are
    obtained from Pfaffians of matrices of the form $S$ and
    \[
    \bigp{
      \begin{array}{ccc}
        T & b & ya \\
        -b^T & 0 & y \\
        -ya^T & -y & 0
      \end{array}
    }.
    \]
    respectively.  Now, since Pfaffians are unchanged by adding a
    multiple of row $i$ to row $j$ and the same multiple of column $i$
    to column $j$, we obtain
    \[
    \Pf
    \bigp{
      \begin{array}{ccc}
        T & b & ya \\
        -b^T & 0 & y \\
        -ya^T & -y & 0
      \end{array}
    }
    =
    \Pf
    \bigp{
      \begin{array}{ccc}
        S & b & 0 \\
        -b^T & 0 & y \\
        0 & -y & 0
      \end{array}
    }
    \]
    which, upon expansion, is $y\Pf(S)$.  So the sign is multiplied by
    the sign of $x$ (which is also the sign of $y$).
  \item The basis contains $n$ but not $n+1$.  Thus the signs are
    given by Pfaffians of matrices of the forms
    \[
    \bigp{
      \begin{array}{cc}
        S & a \\
        -a^T & 0 \\
      \end{array}
    } \quad\hbox{and}\quad
    \bigp{
      \begin{array}{cc}
        T & ya \\
        -ya^T & 0 \\
      \end{array}
    }.
    \]
    Again, using row operations, we obtain
    \[
    \Pf
    \bigp{
      \begin{array}{cc}
        T & ya \\
        -ya^T & 0 \\
      \end{array}
    }
    =
    \Pf
    \bigp{
      \begin{array}{cc}
        S & ya \\
        -ya^T & 0 \\
      \end{array}
    },
    \]
    and so again the sign is multiplied by the sign of $y$ (which is
    also the sign of $x$).
  \item The basis contains $n+1$ but not $n$; this is similar to case
    2.
  \item The basis contains both $n$ and $n+1$.  This is similar to
    case 1.
  \end{enumerate}
  This completes the proof.
\end{proof}

The above results mean that, since $B_n$-represented matroids correspond to
$D_n$-represented matroids of one dimension larger with essentially
the same bases, there is nothing new to be gained by studying their
orientations.

It is natural to wonder, given that our $B_n$-represented matroid is
built from a Lagrangian pair of subspaces, to what extent the signs of
the oriented orthogonal Lagrangian matroids thus represented are
preserved.

\begin{thm}
  The signs of the Lagrangian pair of\/ {\rm (}represented, and so oriented\/{\rm )}
  matroids constituting a $B_n$-represented matroid are the same as
  the signs of the corresponding bases in the $B_n$-representation, up
  to changing sign throughout either constituent.
\end{thm}

\begin{proof}
  Again, we assume that $[n]^*$ and $[n-1]^*\cup n$ are bases.  Thus
  the orthogonal pair are represented by subspaces
  \[
  \bigp{
    \begin{array}{cccc}
      S+ba^T & a & I_{n-1} & -b\\
      -a^T & 0 & 0 & 1
    \end{array}
  }
  \]
  and
  \[
  \bigp{
    \begin{array}{cccc}
      S+ba^T & a & I_{n-1} & -b \\
      b^T & 1 & 0 & 0
    \end{array}
  },
  \]
  using the same notation as in the previous proof.  The signs of these
  matroids are Pfaffian minors of the matrices
  \[
  \bigp{
    \begin{array}{cc}
      S & a \\
      -a^T & 0 \\
    \end{array}
  }
  \quad\hbox{and}\quad
  \bigp{
    \begin{array}{cc}
      T & b \\
      -b^T & 0 \\
    \end{array}
  }
  \]
  indexed by $[n]$ and $[n-1]\cup n^*$ respectively.  Attaching the
  two together to form a $B_n$-representation, we get exactly the
  representation in canonical form from the previous proof, and the proof
  that the signs match up is similar to the calculations there also.
\end{proof}

\vfill
\small
\sc

\noindent
Richard F. Booth\\
Department of Mathematics, UMIST, PO Box 88,Manchester M60~1QD,
United Kingdom\\
\texttt{richard.booth@umist.ac.uk}\\

\medskip
\noindent
Alexandre V. Borovik\\
{\small Department of Mathematics,}
{\small         UMIST, PO Box 88,}
{\small         Manchester M60~1QD,}
{\small         United Kingdom}\\
{\small         \texttt{borovik@umist.ac.uk}}\\

\medskip
\noindent
 Neil White
  {\small Department of Mathematics}
  {\small University of Florida}
  {\small Gainesville, }
  {\small Florida 32611, USA}\\
  {\small \texttt{white@math.ufl.edu}}

\end{document}